\documentclass{amsart}
\usepackage{graphicx}
\usepackage{amssymb}
\newtheorem{thm}{Theorem}[section]

\theoremstyle{definition}

\theoremstyle{question}
\newtheorem{que}[thm]{Question}
\theoremstyle{Conjecture}

\numberwithin{equation}{section}
\begin{document}

\title[NON-NILPOTENT SUBGROUPS OF LOCALLY GRADED
GROUPS]{NON-NILPOTENT SUBGROUPS OF LOCALLY GRADED GROUPS}%
\author{Mohammad Zarrin}%

\address{Department of Mathematics, University of Kurdistan, P.O. Box: 416, Sanandaj, Iran}%
\email{M.Zarrin@uok.ac.ir, Zarrin@ipm.ir}
%\date{}%
%\dedicatory{{\rm }}
%\commby{}%
% ----------------------------------------------------------------
\begin{abstract}

 In this paper, we show that a locally graded group with a finite number $m$ of non-(nilpotent of
class at most $n$) subgroups is (soluble of class at most
$[\log_2(n)]+m+3$)-by-(finite of order $\leq m!$). Also we show
that the derived length of a soluble group with a finite number
$m$ of non-(nilpotent of class at most $n$) subgroups, is at most $[\log_2(n)]+m+1$.
\end{abstract}
\maketitle
% ----------------------------------------------------------------

\section{\textbf{ Introduction and results}}
Let $G$ be a group. A non-nilpotent finite group whose proper
subgroups are all nilpotent is well-known (called Schmidt group).
In 1924, O.Yu. Schmidt studied such groups and proved that such
groups are soluble \cite{Sch}. Subsequently, Newman and Wiegold in
\cite{New}, discussed infinite non-nilpotent groups whose proper
subgroups are all nilpotent. Such groups need not be soluble in
general. For example, the  $Tarski~ Monsters$, which are infinite
simple groups with all proper subgroups of a fixed prime order.

Following \cite{zar} we say that a group $G$ is a
$\mathcal{S}^m$-group if $G$ has exactly $m$ non-nilpotent
subgroups.  More recently Zarrin in \cite{zar} generalized
Schmidt's Theorem and proved that every finite
$\mathcal{S}^m$-group with $m<22$ is soluble. Let $n$ be a
non-negative integer. We say that a group $G$ is an
$\mathcal{S}^m_n$-group, if $G$ has exactly $m$ non-(nilpotent of
class at most $n$) subgroups. Clearly, every
$\mathcal{S}^m_n$-group is a $\mathcal{S}^{r}$-group, for some
$r\leq m$. Here, we show that every locally graded
  group with a finite number $m$ of
non-(nilpotent of class at most $n$) subgroups, is
soluble-by-finite.  Recall that a group $G$ is locally graded if
every non-trivial finitely generated subgroup of $G$ has a
non-trivial finite homomorphic image. This is a rather large class
of groups, containing for instance all residually finite groups
and all locally–(soluble-by-finite) groups.\\

\noindent{\bf Theorem A.}
 Every locally graded $\mathcal{S}^m_n$-group is (soluble of class
at most $[\log_2(n)]+m+3$)-by-(finite of order $\leq m!$).\\

This result suggests that the behavior of non-(nilpotent of class
at most $n$) subgroups has a strong influence on the structure of
the group.\\

 Finding a upper bound for the solubility length of a soluble
group is an important problem in the theory of groups, for example
see \cite{zar2}. It is well-known that a nilpotent group of class
$n$ (or a group without non-(nilpotent of class at most $n$)
subgroups) has derived length  $\leq [\log_2(n)]+1$ (see \cite{d.j.r}, Theorem 5.1.12). Here, we
obtain a result which is of independent interest, namely, the
derived length of  soluble $\mathcal{S}^m_n$-groups is bounded in
terms of $m$ and $n$. (Note that every nilpotent group of class
$n$ is a $\mathcal{S}^m_n$-group
with $m=0$.)\\

 \noindent{\bf Theorem B.}
 Let $G$ be a  soluble $\mathcal{S}^m_n$-group  and $d$ be the
derived length of $G$ . Then $d \leq [\log_2(n)]+m+1.$\\

\section{\textbf{Proofs}}

If $G$ is an arbitrary group, the $norm~B_1(G)$ of $G$ is the
intersection of the normalizers of all subgroups of $G$ and $W(G)$
is the intersection of the normalizers of all subnormal subgroups
of $G$. In 1934 and 1958, respectively, those  concepts were considered
by R. Baer and Wielandt (see also \cite{BHN, BR, LS}). More
recently Zarrin generalized this concept in \cite{zar1}. Here we
define $A_n(G)$ as the intersection of all the normalizers of
non-(nilpotent of class at most $n$) subgroups of $G$, i.e.,
$$A_n(G)=\bigcap_{H\in \mathfrak{T}_n(G)}N_G(H),$$
where  $\mathfrak{T}_n(G)=\big\{H \mid H \text{~is~ a
non-(nilpotent of class at most n) subgroup~of~} G\}$ (with the
stipulation that $A_n(G)= G$ if all subgroups of $G$ are nilpotent
of class at most $n$). Clearly $$B_1(G)\leq A_i(G)\leq
A_{i+1}(G).$$ Moreover, in view of the proof of Theorem A, below,
we can see that, for every locally graded group $G$, we have
$$ A_n(G)\text{ is a soluble normal subgroup of} ~~G~~ \text{of class} \leq [\log_2(n)]+4.$$
  \\

\noindent{\bf{Proof of Theorem A.}}  The group $G$ acts on the
set $\mathfrak{T}_n(G)$
 by conjugation. By
assumption $|\mathfrak{T}_n(G)|=m$. It is easy to see that the
subgroup $A_n(G)$ is the kernel of this action and so $A_n(G)$ is
normal in $G$ and $G/A_n(G)$ is embedded in $ S_{m}$, the
symmetric group of degree $m$. So $$|G/A_n(G)|\leq m!.$$

 Therefore
to complete the proof it is enough to show that $H=A_n(G)$ is
soluble of class at most $[\log_2(n)]+4$. To see this, it is
enough to show that $K=H^{(3)}$ is nilpotent of class at most $n$.
 Suppose on the contrary that $K$ is not nilpotent of
class at most $n$. It follows that every subgroup containing $K$
is not nilpotent of class at most $n$ and so, by definition of
$A_n(G)$, it is a normal subgroup of $H$. Therefore every subgroup of
$H/K$ is normal. That is, $H/K$ is a Dedekind group and so, it
is well-known (see \cite{d.j.r}, Theorem 5.3.7), that $H/K$ is metabelian. From which it follows that
$$H^{(2)}=H^{(3)}=K.\eqno(\bullet)$$ We claim the following conclusions.\\\\
$\mathbf{Step 1.}$ Every proper normal subgroup of $K$ is
nilpotent of class at most $n$.\\

 Suppose, a contrary, that there
exists a proper normal subgroup $M$ of $K=H^{(2)}$ such that $M$
is not nilpotent of class at most $n$. Then we can obtain, by
definition of $A_n(G)$, that $H^{(2)}/M$ is a Dedekind group (so
it is metabelian) and so, in view of $(\bullet)$,
$H^{(2)}=M$, a contradiction.\\\\
$\mathbf{Step 2.}$ The product of all proper normal subgroups of
$K$, say $R$, is a proper nilpotent subgroup of $K$ of
class at most $n$.\\

Suppose that $M_1, M_2,\ldots, M_t$ are proper normal subgroups of
$H^{(2)}$. Then, by step 1, every $M_i$ is soluble and so
$M_1M_2\ldots M_t$ is soluble. Now by $(\bullet)$, we conclude
that $H^{(2)}\neq M_1M_2\ldots M_t$. Therefore $M_1M_2\ldots M_t$
is a proper normal subgroup of $H^{(2)}$ and so, by step 1, it is
nilpotent of class at most $n$. Therefore $R$ is a locally
nilpotent of class at most n group, and so $R$ is nilpotent of
class at most $n$ (note that the class of nilpotent groups of
class at most $n$ is locally closed).
 Also as $(\bullet)$, we have $R\neq H^{(2)}$.\\\\
$\mathbf{Step 3.}$ Finishing the proof. \\

 We note that, by definition of $A_n(G)$, every subgroup of $H^{(2)}$ which is not nilpotent of class at most $c$  is a normal
  subgroup of $H^{(2)}$.
 It follows, as $H^{(2)}/R$ is a simple group, that all proper subgroups of $H^{(2)}/R$ are nilpotent of class at most
 $n$. Since $H^{(2)}$ is locally graded, by the main result of \cite{LMS}, $H^{(2)}/R$ is locally graded. Therefore
 if $H^{(2)}/R$ is
 finitely generated then it must be finite. Thus, by Schmidt's
 Theorem, $H^{(2)}/R$ is soluble, which is contrary to
 $(\bullet)$. If $H^{(2)}/R$ is not finitely generated, then $H^{(2)}/R$ is locally nilpotent of class at most $n$ and so
 $H^{(2)}/R$ is nilpotent of class at most $n$, a contradiction.\\

Now we prove Theorem B.\\

\noindent{\bf{Proof of Theorem B.}}
 Assume that a soluble
group $G$ has derived length $> [\log_2 n]+1+m$ for some $n, m\ge
1$. Then obviously the $m+1$ derived subgroups $G, G',\dots,
G^{(m)}$ are all pairwise distinct and have solubility length $>
[\log_2 n]+1$. Therefore they cannot be nilpotent of class at most
$n$. This shows that $G$ cannot be a $S_n^m$-group, a
contradiction.\\

 Finally, as every
$\mathcal{S}^m_n$-group is a $\mathcal{S}^{r}$-group, for some
$r\leq m$, and by the main result in \cite{zar}, we can see that every
$\mathcal{S}^m_n$-group with $m\leq 21$ is soluble. Hence the
following question arises naturally:

\begin{que}
Assume that $G$ is a $\mathcal{S}^m_n$-group. What relations
between $m, n$ guarantee that $G$ is soluble?
\end{que}

{\noindent\bf Acknowledgements.} I would like to thank the referee for his/her helpful comments.

\end{document}